\documentclass[12pt,dvipdfmx]{article}
\usepackage{tikz}
\usetikzlibrary{intersections,calc,arrows}
\usepackage{amssymb}
\date{}
\makeatletter
\newcommand{\figcaption}[1]{\def\@captype{figure}\caption{#1}}
\newcommand{\tblcaption}[1]{\def\@captype{table}\caption{#1}}

\@addtoreset{equation}{section}
\makeatother
\newcommand{\qed}{\hbox{\rule[-2pt]{3pt}{6pt}}}

\begin{document}
\title {\bf Asymptotics of solution curves of Kirchhoff type elliptic equations 
with logarithmic Kirchhoff function}

\author{Tetsutaro Shibata 
\\
Laboratory of Mathematics, 
\\
Graduate School of Advanced Science and Engineering
\\
Hiroshima University, 
Higashi-Hiroshima, 739-8527, Japan
}

\maketitle
\footnote[0]{E-mail: tshibata@hiroshima-u.ac.jp}
\footnote[0]{This work was supported by JSPS KAKENHI Grant Number JP21K03310.}

\begin{abstract}
We study the one-dimensional nonlocal elliptic equation
\begin{eqnarray*}
-\log(a\Vert u'\Vert_2^2 + b\Vert u\Vert_2^2 + 1)u''(x)&=& \lambda u(x)^p , 
\enskip x \in I:= (0,1), \enskip u(x) > 0, 
\enskip x\in I, 
\\
u(0) &=& u(1) = 0,
\end{eqnarray*}
where $a \ge 0, b > 0, p > 1$ are given constants and $\lambda > 0$ is a bifurcation parameter.  
We establish the precise asymptotic formulas for $u_\lambda(x)$ as $\lambda \to \infty$.

\end{abstract}



\section{Introduction} 		      

We consider the following one-dimensional nonlocal elliptic equation

\begin{equation}
\left\{
\begin{array}{l}
-\log(a\Vert u'\Vert_2^2 + b\Vert u\Vert_2^2 + 1) u''(x)= \lambda u(x)^p, \enskip 
x \in I:= (0,1),
\vspace{0.1cm}
\\
u(x) > 0, \enskip x\in I, 
\vspace{0.1cm}
\\
u(0) = u(1) = 0,
\end{array}
\right.
\end{equation}
where $a \ge 0, b >0, p > 1$ are given constants, $\lambda > 0$ is a bifurcation parameter and 
$\Vert \cdot\Vert_2$ denotes the usual $L^2$-norm.    
\noindent
Equation (1.1) is the nonlocal problem of Kirchhoff type, which is motivated by the following problem (1.2) in [7]:
\begin{equation}
\left\{
\begin{array}{l}
-A\left(\displaystyle{\int_0^1} \vert u(x)\vert^q dx\right)u''(x)= 
\lambda f(x, u(x)), \enskip 
x \in I,
\vspace{0.1cm}
\\
u(x) > 0, \enskip x\in I, 
\vspace{0.1cm}
\\
u(0) = u(1) = 0,
\end{array}
\right.
\end{equation}
where $A = A(w) \ge 0$, which is called Kirchhoff function, 
is a continuous function of $w \ge 0$,

Nonlocal problems have been studied by many investigators, since many problems come from the phenomena of, for instance, biological problems such as population dynamics. Moreover, nonlocal problems have been also derived from numerous physical models and the other area of science, and 
have been studied intensively. We refer to [1--4, 6--11, 13--19], and the references therein. It seems that 
the main interests in this area are existence, nonexistence and the number of positive solutions. 
On the other hand, as far as the author knows, there are a few works which treat (1.1) as the  
bifurcation problems. We refer to [14, 19] and the references therein. In [19], 
the case where $A(\Vert u'\Vert_q) = a\Vert u' \Vert_2^2 + b$ and $f(x, u) = u^p$ in (1.2) has been considered and the existence of a branch of 
positive solutions bifurcating from infinity at $\lambda = 0$ was studied. 
In this paper, we concentrate on the generalized lonlocal Emden--Fowler equation with logarithmic Kirchhoff function, and establish the 
precise global structure of solution curves $u_\lambda(x)$ as $\lambda \to 0$ and 
$\lambda \to \infty$. 
It should be mentioned that, in many cases, Kirchhoff function $A$ contains only 
one of $\Vert u'\Vert_2$ or 
$\Vert u\Vert_2$. However, the Kirchhoff function in (1.1) contains both of them simultaneously. 
It seems that there are few papers which treat such problems as (1.1). Therefore, little is known 
about the properties of the solutions of (1.1).

To state our results, we prepare the following notation. For $p> 1$, let 
\begin{equation}
\left\{
\begin{array}{l}
-W''(x)= W(x)^p, \enskip 
x \in I,
\vspace{0.1cm}
\\
W(x) > 0, \enskip x\in I, 
\vspace{0.1cm}
\\
W(0) = W(1) = 0.
\end{array}
\right.
\end{equation}
We know from [6] that there exists a unique solution $W_p(x)$ of (1.3). 
For $m \ge 1$ and $q \ge 0$, we put 
 \begin{eqnarray}
L_{p,q}:= \int_0^1 \frac{s^q}{\sqrt{1-s^{p+1}}}ds, \enskip
 M_{p,m}:= \int_0^1 (1-s^{p+1})^{(m-1)/2}ds.
 \end{eqnarray}
 Then we have 
 \begin{eqnarray}
\Vert W_p'\Vert_m^m &=& 2^{mp/(p-1)}(p+1)^{m/(p-1)}L_{p,0}^{(mp+m-p+1)/(p-1)}M_{p,m},
\\
\Vert W_p\Vert_\infty &=&  (2(p+1))^{1/(p-1)}L_{p,0}^{2/(p-1)}.
\end{eqnarray} 
 (1.5) and (1.6) have been given in [16]. For completeness, the proof of (1.5) and (1.6) will be given in Appendix. 

\vspace{0.2cm}

In the following Theorems 1.1--1.3, we consider the case $a = 0$ and $b$, which is rewritten by $d$ in (1.1), since it is convenient to distinguish between the constant in front of  $\Vert u\Vert_2$ in the 
case $a = 0$ and $a \not=0$. Namely, we first consider the equation 
\begin{equation}
\left\{
\begin{array}{l}
-\log(d\Vert u\Vert_2^2 + 1) u''(x)= \lambda u(x)^p, \enskip 
x \in I:= (0,1),
\vspace{0.1cm}
\\
u(x) > 0, \enskip x\in I.
\vspace{0.1cm}
\\
u(0) = u(1) = 0,
\end{array}
\right.
\end{equation}
where $d > 0$ is a given constant. Now we state our main results. 

\vspace{0.2cm}

\noindent
{\bf Theorem 1.1.} {\it  Consider (1.7). Assume that $p > 3$. 
Then for any $\lambda > 0$, there exists a unique solution $u_\lambda$ of (1.7). Further, 
as $\lambda \to \infty$,
\begin{eqnarray}
u_\lambda(x) &=& \left(\frac{\lambda\Vert W_p\Vert_2^{1-p}}{d}\right)^{1/(3-p)}
\left(1 + \frac{1}{2(3-p)}d^{(p-1)/(p-3)}(\lambda
\Vert W_p\Vert_2^{1-p})^{2/(3-p)}(1 + o(1))\right)
\nonumber
\\
&&\times\Vert W_p\Vert_2^{-1}W_p(x).
\end{eqnarray}
Moreover, as $\lambda \to 0$,
\begin{eqnarray}
u_\lambda(x) &=& \left(\frac{2}{p-1}\right)^{1/(p-1)}
\lambda^{-1/(p-1)}\left(\log\frac{1}{\lambda}\right)^{1/(p-1)}
\\
&& \qquad \qquad \qquad \times
\left\{1 + \frac{1}{p-1}\frac{\log\log\frac{1}{\lambda}}{\log\frac{1}{\lambda}}
(1 + o(1))\right\}W_p(x).
\nonumber
\end{eqnarray}
}

\vspace{0.2cm}

\noindent
{\bf Theorem 1.2.} {\it  Consider (1.7). Assume that $1 < p < 3$. Furthermore, put
\begin{eqnarray}
\nu:= \displaystyle{\frac{2}{p-1}\frac{dt_2^{(3-p)/2}}{dt_2 + 1}}\Vert W_p\Vert_2^{p-1},
\end{eqnarray}
where $t_2 > 0$ is a constant deretmined laler. 
Then the following three cases occur:

\noindent
(i) If $0 < \lambda < \nu$, then there exist 
exactly two solutions $u_{\lambda,1}$ and $u_{\lambda,3}$ of (1.7) satisfying 
$u_{\lambda,1} < u_{\lambda,3}$ in $I$. 
Furthermore, as $\lambda \to 0$, 
\begin{eqnarray}
u_{\lambda, 1}(x) &=& \left(\frac{\lambda \Vert W_p\Vert_2^{1-p}}{d}\right)^{1/(3-p)}
 \left(1 + \frac{1}{2(3-p)}\left(\lambda \Vert W_p\Vert_2^{1-p}\right)^{2/(3-p)}d^{(1-p)/(3-p)}
 (1 + o(1))\right)
 \nonumber
 \\
 &&\qquad \qquad \qquad \qquad \qquad \times \Vert W_p\Vert_2^{-1}W_p(x),
\\
u_{\lambda,3}(x) &=& 
\left(\frac{2}{p-1}\right)^{1/(p-1)}\lambda^{-1/(p-1)}\left(
 \log\frac{1}{\lambda}\right)^{1/(p-1)}
 \\
 &&\qquad \qquad \qquad \qquad \times
 \left\{1 + \frac{1}{p-1}\frac{\log(\log\frac{1}{\lambda})}{ \log\frac{1}{\lambda}}(1 + o(1))
 \right\}W_p(x).
 \nonumber
\end{eqnarray}

\noindent
(ii) If $\lambda = \nu$, then there exists exactly one solution $u_\lambda(x)$ 
of (1.7). 

\noindent
(iii) If $\lambda > \nu$, then there exists no solution $u_\lambda(x)$ 
of (1.7). 
}
 
\vspace{0.2cm} 
 
 \noindent
 {\bf Theorem 1.3.} {\it Consider (1.7). Let $p=3$. 

\noindent
(i)  Let 
$\lambda\ge d\Vert W_3\Vert_2^{2} $. Then there exists no solution of (1.7). 

\noindent
(ii)  Let 
$0 < \lambda < d\Vert W_3\Vert_2^{2}$. Then there exists exactly one solution $u_\lambda$ 
of (1.7). 

\noindent
(iii) 
Let $0 < \lambda < d\Vert W_3\Vert_2^{2}$. Then as $\lambda \to d\Vert W_3\Vert_2^{2}$,
\begin{eqnarray}
u_\lambda(x) &=& 
\frac{\sqrt{2}}{d}
 \sqrt{d- \lambda\Vert W_3\Vert_2^{-2} }\left\{1 + \frac{2}{3}
 (d- \lambda\Vert W_3\Vert_2^{-2})(1 + o(1))
 \right\}
 \Vert W_3\Vert_2^{-1}W_3(x).
\end{eqnarray}
Furthermore, as $\lambda \to 0$,
\begin{eqnarray}
u_\lambda(x) 
&=& \lambda^{-1/2}\left(
 \log\frac{1}{\lambda}\right)^{1/2}
 \left\{1 +\frac{\log(\log\frac{1}{\lambda})}{ 2\log\frac{1}{\lambda}}(1 + o(1))
 \right\}W_p(x).
\end{eqnarray}
 }
  
\vspace{0.2cm}

Now we consider the equation (1.1). We show that (1.1) is reduced to (1.7).  

\vspace{0.2cm}

\noindent
{\bf Theorem 1.4.} {\it Consider (1.1). Then (1.1) is reduced to (1.7)  with $d = d_0$, where 
\begin{eqnarray}
d_0 := \frac{4aL_{p,0}^2M_{p,2}}{L_{p,2}} + b.
\end{eqnarray} 
Namely, the solution $u_\lambda$ of (1.1) with $a > 0$ and $b >0$ satisfies
(1.7) with $d = d_0$. Therefore, the solution $u_\lambda$ of (1.1) satisfies all the results 
in Theorems 1.1--1.3 
with $d_0$. }
  
\vspace{0.2cm}

The rest of this paper is organized as follows. In Sec. 2, we prove Theorems 1.1--1.3 
by using the argument in [1] and time map method (cf. [12]). In Sect. 3, we prove Theorem 1.4. 
The final section is the Appendix, in which the proofs of (1.5) and (1.6) will be given for the reader's 
convenience.

\section{Proofs of Theorems 1.1--1.3}

In this section, we consider (1.7). 
By [5], we know that if $u_\lambda$ is a solution of (1.7), then $u_\lambda$ satisfies 
\begin{eqnarray}
u_\lambda(x) &=& u_\lambda(1-x), \quad 0 \le x\le \frac12,
\\
\alpha&:=& \Vert u_\lambda\Vert_\infty = u_\lambda\left(\frac12\right),
\\
u_\lambda'(x) &>& 0, \quad 0 \le x < \frac12.
\end{eqnarray}
For a given $\lambda > 0$, 
let $w_\lambda(x)$ be a unique solution of 
\begin{equation}
\left\{
\begin{array}{l}
-w''(x)= \lambda w(x)^p, \enskip 
x \in I,
\vspace{0.1cm}
\\
w(x) > 0, \enskip x\in I, 
\vspace{0.1cm}
\\
w(0) = w(1) = 0.
\end{array}
\right.
\end{equation}
It is clear that  $w_\lambda = \lambda^{-1/(p-1)}W_p$. We explain the existence of the solutions $u_\lambda$ of (1.7) by using the idea in [1]. We put $M(t):= \log(dt + 1)$ and 
consider the equation for $t > 0$:
\begin{eqnarray}
M(t) = \log(dt + 1) = \Vert w_\lambda\Vert_2^{1-p}t^{(p-1)/2}.
\end{eqnarray}
Assume that $t_\lambda > 0$ satisfies (2.5). We put 
$\gamma:= t_\lambda^{1/2}\Vert w_\lambda\Vert_2^{-1}$ 
and 
\begin{eqnarray}
u_\lambda := \gamma w_\lambda = t_\lambda^{1/2}\Vert w_\lambda\Vert_2^{-1}w_\lambda
=  t_\lambda^{1/2}\Vert W_\lambda\Vert_2^{-1}W_\lambda.
\end{eqnarray}
Then by (2.5),  we have 
$M(\Vert \gamma w_\lambda\Vert_2^2) = M(t_\lambda) = \gamma^{p-1}$. 
Then we have 
\begin{eqnarray}
-M(\Vert u_\lambda\Vert_2^2)u_\lambda''(x) &=& -M(\Vert\gamma  w_\lambda\Vert_2^2)\gamma w_\lambda''(x)
\\
&=& \gamma^p\lambda w_\lambda^p = \lambda(\gamma w_\lambda(x))^p 
\nonumber
\\
&=& \lambda u_\lambda(x)^p.
\nonumber
\end{eqnarray}
Let 
\begin{eqnarray} 
f(t):= \frac{\log(dt + 1)}{t^{(p-1)/2}}. 
\end{eqnarray}
By (2.5) and (2.7), to find the solutions of (1.7), we look for solutions $t_\lambda > 0$  of the following equation of $t > 0$:
\begin{eqnarray}
f(t) = \lambda\Vert W_p\Vert_2^{1-p}. 
\end{eqnarray}
Assume that $u_\lambda$ is a solution of (1.7). Then 
$u_\lambda$ is a solution of (2.4) with $\lambda/M(\Vert u_\lambda\Vert_2^2)$. 
Therefore, by the uniqueness of $W_p$ in (1.3), there exists a unique constant $\Lambda > 0$ such that $u_\lambda = \Lambda w_\lambda$. Then we see that 
$\Lambda = \Vert u_\lambda\Vert_2\Vert w_\lambda\Vert_2^{-1}$. Then 
we put $u_\lambda = \Vert u_\lambda\Vert_2\Vert w_\lambda\Vert_2^{-1}w_\lambda$ and 
$t_\lambda :=\Vert u_\lambda\Vert_2^2$. Since $u_\lambda$ satisfies (1.7), we have 
\begin{eqnarray}
-M(\Vert u_\lambda\Vert_2^2)\Vert u_\lambda\Vert_2\Vert w_\lambda\Vert_2^{-1}w_\lambda''(x) &=& \lambda \Vert u_\lambda\Vert_2^p\Vert w_\lambda\Vert_2^{-p}w_\lambda(x)^p.
\end{eqnarray}
This implies that 
\begin{eqnarray}
M(t_\lambda) = M(\Vert u_\lambda\Vert_2^2) = 
\Vert u_\lambda\Vert_2^{p-1}\Vert w_\lambda\Vert_2^{1-p} = \lambda\Vert W_p\Vert_2^{1-p}\Vert u_\lambda\Vert_2^{p-1} = \lambda\Vert W_p\Vert_2^{1-p}t_\lambda^{(p-1)/2}.
\nonumber
\end{eqnarray}
This implies (2.5). Therefore, the solution of (1.7) coincides with the solution $t$ of (2.5). 
\vspace{0.2cm}

\noindent
{\bf Lemma 2.1.} {\it Assume that $p > 3$. Then for any given $\lambda > 0$, 
there exists a unique $t_\lambda > 0$ such that 
$f(t_\lambda) = \lambda\Vert W_p\Vert_2^{1-p}$. 
Furthermore, as $\lambda \to \infty$,
\begin{eqnarray}
u_\lambda(x) &=& \left(\frac{\lambda\Vert W_p\Vert_2^{1-p}}{d}\right)^{1/(3-p)}
\left(1 + \frac{1}{2(3-p)}d^{(p-1)/(p-3)}(\lambda
\Vert W_p\Vert_2^{1-p})^{2/(3-p)}(1 + o(1))\right)
\nonumber
\\
&&\times\Vert W_p\Vert_2^{-1}W_p(x).
\end{eqnarray}
Moreover, as $\lambda \to 0$,
\begin{eqnarray}
u_\lambda(x) &=& 
\left(\frac{2}{p-1}\right)^{1/(p-1)}
\lambda^{-1/(p-1)}\left(\log\frac{1}{\lambda}\right)^{1/(p-1)}
\\
&&\qquad \qquad \qquad \times 
\left\{1+ \frac{1}{p-1}\frac{\log\left(\log\frac{1}{\lambda}\right)}{\log\left(\frac{1}{\lambda}\right)}(1 + o(1)) \right\}W_p(x).
\nonumber
\end{eqnarray}
}
\noindent
{\bf Proof.} By (2.8), we have 
\begin{eqnarray}
f'(t) &=& -\frac{p-1}{2}t^{-(1+p)/2}\log(dt + 1) + t^{(1-p)/2}\frac{d}{dt+1}
\\
&=& t^{-(1+p)/2}g(t),
\nonumber
\end{eqnarray}
where 
\begin{eqnarray}
g(t):= -\frac{p-1}{2}\log(dt + 1) + \frac{dt}{dt+1}.
\end{eqnarray}
By direct calculation, we have 
\begin{eqnarray}
g'(t) = \frac{d}{(dt+1)^2}\left\{-\frac{p-1}{2}(dt+1) + 1\right\}.
\end{eqnarray}
By this, we see that $g'(t) < 0$ for $t > 0$. We know that $g(0) = 0$. 
Therefore, $g(t) < 0$ for $t > 0$. By this and (2.13), we see that $f'(t) < 0$ for $t > 0$. Namely, 
$f(t)$ is strictly decreasing for $t > 0$. Further, 
$\lim_{t \to 0}f(t) = \infty$ and $\lim_{t \to \infty}f(t) = 0$. 
Therefore, there exists a unique $t_\lambda > 0$ such that the equation (2.9) holds.

\noindent
We first assume that $\lambda \to \infty$. We see that $t_\lambda \to 0$ as $\lambda \to \infty$. By this 
and Taylor expansion, we have 
\begin{eqnarray}
\log(dt_\lambda + 1) = dt_\lambda -\frac12d^2t_\lambda^2 (1 + o(1)).
\end{eqnarray} 
By this and (2.9), we have 
\begin{eqnarray}
f(t_\lambda) = \frac{\log(dt_\lambda + 1)}{t_\lambda^{(p-1)/2}} 
 &=& dt_\lambda^{(3-p)/2} - \frac12d^2t_\lambda^{(5-p)/2}(1 + o(1)) = 
 \lambda\Vert W_p\Vert_2^{1-p}. 
\end{eqnarray}  
We put 
\begin{eqnarray}
t_\lambda =\left(\frac{\lambda\Vert W_p\Vert_2^{1-p}}{d}\right)^{2/(3-p)}(1 + R_\lambda),
\end{eqnarray} 
where $R_\lambda$ is the remainder term, which satisfies $R_\lambda \to 0$ as 
$\lambda \to \infty$. By this and (2.17), we have 
\begin{eqnarray}
&&\lambda\Vert W_p\Vert_2^{1-p}(1 + R_\lambda)^{(3-p)/2} - \frac12d^{(1-p)/(3-p)}
(\lambda\Vert W_p\Vert_2^{1-p})^{(5-p)/(3-p)}(1 + R_\lambda)^{(5-p)/2} 
\\
&&= \lambda
\Vert W_p\Vert_2^{1-p}.
\nonumber
\end{eqnarray}
Then by (2.19) and 
Taylor expansion, we have 
\begin{eqnarray}
\frac{3-p}{2}\lambda
\Vert W_p\Vert_2^{1-p}R_\lambda = \frac{1}{2}d^{(1-p)/(3-p)}
(\lambda\Vert W_p\Vert_2^{1-p})^{(5-p)/(3-p)}
(1 + o(1)).
\end{eqnarray}
This implies 
\begin{eqnarray}
R_\lambda = \frac{1}{3-p}d^{(p-1)/(p-3)}(\lambda
\Vert W_p\Vert_2^{1-p})^{2/(3-p)}(1 + o(1)).
\end{eqnarray}
By (2.6), (2.18), (2.21) and Taylor expansion,  as $\lambda \to \infty$, we obtain
\begin{eqnarray}
u_\lambda(x) &=& t_\lambda^{1/2}\Vert W_p\Vert_2^{-1}W_p(x)
\\
&=& \left(\frac{\lambda\Vert W_p\Vert_2^{1-p}}{d}\right)^{1/(3-p)}
\left(1 + \frac{1}{2(3-p)}d^{(p-1)/(p-3)}(\lambda
\Vert W_p\Vert_2^{1-p})^{2/(3-p)}(1 + o(1))\right)
\nonumber
\\
&&\times\Vert W_p\Vert_2^{-1}W_p(x).
\nonumber
\end{eqnarray}
This implies (2.11). Next, we assume that $\lambda \to 0$ and show (2.12). By Fig. 3, 
it is clear that $t_\lambda \to \infty$ as $\lambda \to 0$. We look for $t_\lambda$ of the form 
$t_\lambda = C\lambda^{-2/(p-1)}\left(\log\frac{1}{\lambda}\right)^q(1 + R_\lambda)$, where 
$C$ and $q$ are constants and $R_\lambda$ is the remainder term satisfying $R_\lambda \to 0$ as $\lambda \to 0$. By Taylor expansion, we have 
\begin{eqnarray}
\log(dt_\lambda + 1) &=&\log t_\lambda + O(1) 
\\
&=& \frac{2}{p-1}\log\frac{1}{\lambda} + q\log\left(\log\frac{1}{\lambda}\right) + R_\lambda + O(1),
\nonumber
\\
t_\lambda^{(p-1)/2}\lambda\Vert W_p\Vert_2^{1-p} 
&=& \lambda\left\{C\lambda^{-2/(p-1)}\left(\log\frac{1}{\lambda}\right)^q
(1 + R_\lambda)\right\}^{(p-1)/2}
\Vert W_p\Vert_2^{1-p} 
\\
&=& \Vert W_p\Vert_2^{1-p} C^{(p-1)/2}\left(\log\frac{1}{\lambda}\right)^{q(p-1)/2}
\left\{1 + \frac{p-1}{2}R_\lambda + o(R_\lambda)\right\}.
\nonumber
\end{eqnarray}
This implies $C = (2/(p-1))^{2/(p-1)}\Vert W_p\Vert_2^{2}$ and $q = 2/(p-1)$. 
By this, (2.23) and (2.24), we have 
\begin{eqnarray}
\frac{2}{p-1}\log\left(\log\frac{1}{\lambda}\right) + R_\lambda + O(1) = 
R_\lambda(1 + o(1))\log\frac{1}{\lambda}.
\end{eqnarray}
This implies that 
\begin{eqnarray}
R_\lambda = \frac{2}{p-1}\frac{\log\left(\log\frac{1}{\lambda}\right)}{\log\left(\frac{1}{\lambda}\right)}(1 + o(1)).
\end{eqnarray}
By this and (2.6), as $\lambda \to 0$, 
\begin{eqnarray}
u_\lambda(x) &=& t_\lambda^{1/2}\Vert W_p\Vert_2^{-1}W_p(x) 
\\
&=& C^{1/2}\lambda^{-1/(p-1)}\left(\log\frac{1}{\lambda}\right)^{q/2}
\left(1 + R_\lambda + o(R_\lambda)\right)^{1/2}\Vert W_p\Vert_2^{-1}W_p(x) 
\nonumber
\\
&=& \left(\frac{2}{p-1}\right)^{1/(p-1)}
\lambda^{-1/(p-1)}\left(\log\frac{1}{\lambda}\right)^{1/(p-1)}
\left(1+\frac{1}{p-1}\frac{\log\left(\log\frac{1}{\lambda}\right)}{\log\left(\frac{1}{\lambda}\right)}(1 + o(1)) \right)W_p(x).
\nonumber
\end{eqnarray}
This implies (2.12). Thus the proof is complete. \qed

\vspace{0.2cm}

We obtain Theorem 1.1 by Lemma 2.1.  We next prove Theorem 1.2. 

\vspace{0.2cm}

\noindent
{\bf Lemma 2.2.} {\it Assume that $1 < p < 3$. Then there exists a constant 
$\nu$ and the following three cases occur.  

\noindent
(i) if $0 < \lambda < \nu$, then there exist 
two solutions $u_{\lambda,1}$ and $u_{\lambda,3}$ of (1.7). 

\noindent
(ii) If $\lambda = \nu$, then there exists exactly one solution $u_\lambda(x)$ 
of (1.7). 

\noindent
(iii) If $\lambda > \nu$, then there exists no solution $u_\lambda(x)$ 
of (1.7). 
  }
  
\noindent
{\bf Proof.} By (2.15), we see that if $t_0 = (3-p)/(d(p-1))$, then $g'(t_0) = 0$. This implies that 
$g(t)$ is increasing in $0 < t < t_0$ and attains the maximum at $t = t_0$ and decreasing in 
$t > t_0$. Since $g(0) = 0, g(t_0) > 0$ and $g(t) \to -\infty$ as $t \to \infty$, we see that 
there exists $t =t_2$ such that $g(t_2) = 0$. Namely, $f'(t_2) = 0$. We know that 
$\lim_{x \to 0}f(x) = 0$ and $\lim_{x \to \infty}f(x) =0$. Furthermore, $f(t)$ attains the maximum that $t = t_2$ (cf. Fig. 4 below). 
By (2.13) and (2.14), we see that $t_2$ satisfies
\begin{eqnarray}
\log(dt_2 + 1) = \frac{2}{p-1}\frac{dt_2}{dt_2 + 1}.
\end{eqnarray}
We note that there exists exactly one $t_2$ satisfying (2.28). The reason is simple. 
We consider the graph of $h(x):= \log (x+1) - \frac{2x}{(p-1)(x+1)}$. 
Then $h'(x) = \frac{1}{(x+1)^2}(x- \frac{3-p}{p-1})$. Therefore, $h(x)$ is strictly decreasing in 
$0 < x < (3-p)/(p-1)$ and strictly increasing in $x > (3-p)/(p-1)$. 
Further, $h(0) = 0$ and $\lim_{x \to \infty}h(x) = \infty$. Therefore, 
there exists a unique $t_2$ such that $(3-p)/(p-1)< t_2 < C$ satisfying (2.28), 
since (2.28) does not hold for $t_2 \gg 1$. By this, (2.8) and (2.28), we have 
\begin{eqnarray}
f(t_2) = \frac{2}{p-1}\frac{dt_2^{(3-p)/2}}{dt_2 + 1} = \nu\Vert W_p\Vert_2^{1-p},
\end{eqnarray}
where $\nu:= \displaystyle{\frac{2}{p-1}\frac{dt_2^{(3-p)/2}}{dt_2 + 1}}\Vert W_p\Vert_2^{p-1}$. 
If $0 < \lambda < \nu$, then there exist exactly two $t_1 < t_3$ such that $f(t_j) = \lambda\Vert W_p\Vert_2^{p-1}$ ($j = 1,3$).

This implies that the exists exactly two solutions of (1.7) $u_{\lambda,1}$ and $u_{\lambda,3}$ 
corresponding to $t_1$ and $t_3$. Similarly, 
if $\lambda = \nu$, then there exists one solution of (1.7) and if $\lambda > \nu$, then 
there exist no solutions of (1.7). Thus the proof is complete. \qed

\vspace{0.2cm}

Now we consider the asymptotic behavior of $u_{\lambda,1}$ and $u_{\lambda,3}$ as 
$\lambda \to 0$. 

\vspace{0.2cm}

\noindent
{\bf Lemma 2.3.} {\it Assume that $1 < p < 3$. Then as $\lambda \to 0$, 
\begin{eqnarray}
u_{\lambda, 1}(x) &=& \left(\frac{\lambda \Vert W_p\Vert_2^{1-p}}{d}\right)^{1/(3-p)}
 \left(1 + \frac{1}{2(3-p)}\left(\lambda \Vert W_p\Vert_2^{1-p}\right)^{2/(3-p)}d^{(1-p)/(3-p)}
 (1 + o(1))\right)
 \nonumber
 \\
 &&\qquad \qquad \qquad \qquad \qquad \times \Vert W_p\Vert_2^{-1}W_p(x),
\\
u_{\lambda,3}(x) &=& 
\left(\frac{2}{p-1}\right)^{1/(p-1)}\lambda^{-1/(p-1)}
 \log\left(\frac{1}{\lambda}\right)^{1/(p-1)}
 \\
 &&\qquad \qquad \qquad \qquad \times
 \left\{1 + \frac{1}{p-1}\frac{\log(\log\frac{1}{\lambda})}{ \log\frac{1}{\lambda}}(1 + o(1))
 \right\}W_p(x).
 \nonumber
\end{eqnarray}
}
{\bf Proof.} Since $1 < p < 3$, and $t_1 < t_2 < t_3 $, 
we see from Fig. 4 that $t_1 \to 0$ and $t_3 \to \infty$ as $\lambda \to 0$. 
We first prove (2.30). By (2.4) and Taylor expansion, we have 
\begin{eqnarray}
\frac{dt_1 - (dt_1)^2/2 + o(t_1^2)}{t_1^{(p-1)/2}} = 
 dt_1^{(3-p)/2} - \frac12d^2t_1^{(5-p)/2}(1 + o(1)) =  \lambda \Vert W_p\Vert_2^{1-p}.
 \end{eqnarray}
 By this, we have 
 \begin{eqnarray}
 t_1 &=& \left(\frac{\lambda \Vert W_p\Vert_2^{1-p}}{d}\right)^{2/(3-p)}(1 + \eta),
 \end{eqnarray}
 where $\eta \to 0$ as $\lambda \to 0$. By (2.32), (2.33) and Taylor expansion, we have 
 \begin{eqnarray}
&& \lambda \Vert W_p\Vert_2^{1-p}(1 + \eta)^{(3-p)/2} 
 -\frac12d^2\left(\frac{\lambda \Vert W_p\Vert_2^{1-p}}{d}\right)^{(5-p)/(3-p)}
 (1 + o(1)) 
 \\
&& =  \lambda \Vert W_p\Vert_2^{1-p}\left(1 + \frac{3-p}{2}\eta + o(\eta)\right)
  - \frac12d^2\left(\frac{\lambda \Vert W_p\Vert_2^{1-p}}{d}\right)^{(5-p)/(3-p)}
 (1 + o(1)) 
 \nonumber
 \\
&& =  \lambda \Vert W_p\Vert_2^{1-p}.
 \nonumber
 \end{eqnarray}
 This implies that 
 \begin{eqnarray}
 \eta &=& \frac{1}{3-p}\left(\lambda \Vert W_p\Vert_2^{1-p}\right)^{2/(3-p)}d^{(1-p)/(3-p)}
 (1 + o(1)).
 \end{eqnarray}
 By (2.6), (2.33), (2.35) and 
 Taylor expansion, we have
 \begin{eqnarray}
 u_{\lambda,1}(x) &=& t_1^{1/2}\Vert W_p\Vert_2^{-1}W_p(x)
 \\
 &=& \left(\frac{\lambda \Vert W_p\Vert_2^{1-p}}{d}\right)^{1/(3-p)}
 \left(1 + \frac12\eta(1 + o(1))\right)\Vert W_p\Vert_2^{-1}W_p(x)
 \nonumber
 \\
 &=& \left(\frac{\lambda \Vert W_p\Vert_2^{1-p}}{d}\right)^{1/(3-p)}
 \left(1 + \frac{1}{2(3-p)}\left(\lambda \Vert W_p\Vert_2^{1-p}\right)^{2/(3-p)}d^{(1-p)/(3-p)}
 (1 + o(1))\right)
 \nonumber
 \\
 &&\qquad \qquad \qquad \qquad \qquad \times \Vert W_p\Vert_2^{-1}W_p(x).
 \nonumber
 \end{eqnarray} 
 This implies (2.30). We next prove (2.31). Since 
 $t_3 \to \infty$ as $\lambda \to 0$, by (2.8) and (2.9) and Taylor expansion, we have 
 \begin{eqnarray}
\log(bt_3+1) &=& \log t_3 + O(1) = t_3^{(p-1)/2}\lambda \Vert W_p\Vert_2^{1-p}.
 \end{eqnarray}
 Then the situation is the same as that of (2.23), (2.24) and (2.26). Therefore, by the same argument as 
 that to obtain (2.12), we obtain (2.31). So we omit the proof. 
 Thus the proof is complete. \qed

 \vspace{0.2cm}

\noindent
{\bf Lemma 2.4.} {\it Assume that $p=3$. 

\noindent
(i)  Let 
$\lambda \ge d\Vert W_3\Vert_2^{2}$. Then there exists no solution of (1.7). 

\noindent
(ii)  Let 
$0 < \lambda < d\Vert W_3\Vert_2^{2}$. Then there exists exactly one solution $u_\lambda$ 
of (1.7). 

\noindent
(iii) 
Let $0 < \lambda < d\Vert W_3\Vert_2^{2}$. Then as $\lambda \to d\Vert W_3\Vert_2^{2}$,
\begin{eqnarray}
u_\lambda(x) &=& 
\frac{\sqrt{2}}{d}
 \sqrt{d- \lambda\Vert W_3\Vert_2^{-2} }\left\{1 + \frac{2}{3}
 (d- \lambda\Vert W_3\Vert_2^{-2})(1 + o(1))
 \right\}
 \Vert W_3\Vert_2^{-1}W_3(x).
\end{eqnarray}
Furthermore, as $\lambda \to 0$,
\begin{eqnarray}
u_\lambda(x) 
&=& \lambda^{-1/(2)}\left\{
 \log\left(\frac{1}{\lambda}\right)\right\}^{1/2}
 \left\{1 + \frac{\log(\log\frac{1}{\lambda})}{ 2\log\frac{1}{\lambda}}(1 + o(1))
 \right\}W_p(x).
\end{eqnarray}
}

\noindent
{\bf Proof.} We first prove (i) and (ii). By (2.8), we have 
\begin{eqnarray}
f(t) = \frac{\log(dt+1)}{t}.
\end{eqnarray}
Then by (2.14) and (2.15), for $t > 0$, we have 
\begin{eqnarray}
g(t) &=& -\log(dt+1)+ \frac{dt}{dt+1},
\\
g'(t) &=& -\frac{d^2t}{(dt+1)^2} < 0.
\end{eqnarray} 
By this, $g(t)$ is strictly decreasing for $t \ge 0$ and $g(0) = 0$. This implies that $g(t) < 0$ for 
$t > 0$.  By this and (2.13), $f'(t) < 0$ for $t > 0$ and $f(0) = \lim_{t \to 0}f(t) = d > 0$. Further, 
$f(t) \to 0$ as $t \to \infty$. By this and (2.9), we see that if 
$\lambda\Vert W_3\Vert_2^{-2} \ge d$, then there exists no solution $t_\lambda > 0$ of (2.9). 
Further, if $0 < \lambda \Vert W_p\Vert_2^{-2} < d$, then 
there exists exactly one solution $t_\lambda > 0$ of (2.9). Therefore, 
we
obtain (i) and (ii).

We next prove (iii). In this case, it is clear that $t_\lambda \to 0$ as 
$\lambda \to  d\Vert W_3\Vert_2^{-2}$. By this, (2.8), 
(2.9) and Taylor expansion, we have 
\begin{eqnarray}
f(t_\lambda) &=& \frac{\log(dt_\lambda+1)}{t_\lambda} 
=\frac{dt_\lambda - (1/2)d^2t_\lambda^2 + (1/3)d^3t_\lambda^3 + O(t_\lambda^4)}{t_\lambda} 
\\
&=& d - \frac12d^2t_\lambda + \frac13d^3t_\lambda^2 + O(t_\lambda^3) = \lambda \Vert W_3\Vert_2^{-2}.
\nonumber
\end{eqnarray}
 This implies that 
 \begin{eqnarray}
 t_\lambda = \frac{2}{d^2}(d - \lambda\Vert W_3\Vert_2^{-2}) + R,
 \end{eqnarray}
 where $R$ is a remainder term. By this, (2.43) and direct calculation, we have 
 \begin{eqnarray}
 R = \frac{8}{3d^2}(d - \lambda\Vert W_3\Vert_2^{-2})^2(1 + o(1)).
 \end{eqnarray}
 By (2.6), (2.44) and (2.45), we obtain 
 \begin{eqnarray}
 u_\lambda(x) &=& t_\lambda^{1/2}\Vert W_3\Vert_2^{-1}W_3(x) 
 \\
 &=& \frac{\sqrt{2}}{d}
 \sqrt{d- \lambda\Vert W_3\Vert_2^{-2} }\left\{1 + \frac{2}{3}
 (d- \lambda\Vert W_3\Vert_2^{-2})(1 + o(1))
 \right\}
 \Vert W_3\Vert_2^{-1}W_3(x).
 \nonumber
 \end{eqnarray}
 This implies (2.39). 
 Finally, we see that the argument as that to obtain (1.12) is available in the case 
 $p =3$. Therefore, (1.14) follows from (1.12) by putting $p = 3$. Thus the proof is complete. \qed

 \section{Proof of Theorem 1.4}

In this section, we consider (1.1). Let $a > 0$ and $b > 0$ in (1.1). We show that (1.1) is equivalent to (1.7) if we put
\begin{eqnarray}
d = d_0 := a\frac{4aL_{p,0}^2M_{p,2}}{L_{p,2}} + b.
\end{eqnarray}
To do this, we need two lemmas. 

\vspace{0.2cm}

\noindent
{\bf Lemma 3.1.} {\it Assume that $u_\lambda$ satisfies (1.1). Then 
\begin{eqnarray}
\Vert u_\lambda'\Vert_2^2 =  \frac{4L_{p,0}^2M_{p,2}}{L_{p,2}}\Vert u_\lambda\Vert_2^2.
\end{eqnarray}
}
{\bf Proof.} We put $H:= \log(a\Vert u_\lambda'\Vert_2 + b\Vert u_\lambda\Vert_2^2 + 1)$. 
By (1.1), we have 
\begin{eqnarray}
Hu_\lambda''(x) + \lambda u_\lambda(x)^p = 0.
\end{eqnarray}
This implies 
\begin{eqnarray}
\{Hu_\lambda''(x) + \lambda u_\lambda(x)^p\}u_\lambda'(x) = 0.
\end{eqnarray}
We recall that $\alpha = \Vert u_\lambda\Vert_\infty$, which 
is defined in (2.2). By this, (2.2) and putting $x = 1/2$, for $x \in \bar{I}$, we have 
\begin{eqnarray}
\frac12H u_\lambda'(x)^2 + \frac{1}{p+1}\lambda u_\lambda(x)^p = \mbox{constant} 
=  \frac{1}{p+1}\lambda\alpha^{p+1}.
\end{eqnarray}
By this and (2.3), for $0 \le x \le 1/2$, we have 
\begin{eqnarray}
u_\lambda'(x) = \sqrt{k(\alpha^{p+1}-u_\lambda(x)^{p+1})},
\end{eqnarray} 
where $k := 2\lambda/(H(p+1))$. 
By this, (1.4), (2.1), (2.3) and putting $u_\lambda =\theta = \alpha s$, we have 
\begin{eqnarray}
\Vert u_\lambda'\Vert_2^2 &=& 2\int_0^{1/2} \sqrt{k(\alpha^{p+1}-u_\lambda(x)^{p+1})}
u_\lambda'(x)dx \\
&=& 2\int_0^\alpha  \sqrt{k(\alpha^{p+1} - \theta^{p+1})}d\theta
\nonumber
\\
&=& 2\sqrt{k}\alpha^{(p+3)/2}\int_0^1 \sqrt{1-s^{p+1}}ds 
=  2\sqrt{k}M_{p,2}\alpha^{(p+3)/2},
\nonumber
\\
\Vert u_\lambda\Vert_2^2 &=& 2\int_0^{1/2} u_\lambda(x)^2\frac{u_\lambda'(x)}
{\sqrt{k(\alpha^{p+1}-u_\lambda(x)^{p+1})}}dx
\\
&=& 2\int_0^\alpha \frac{\theta^2}{\sqrt{k(\alpha^{p+1}-\theta^{p+1})}}d\theta
\nonumber
\\
&=& \frac{2}{\sqrt{k}}\alpha^{(5-p)/2}\int_0^1 \frac{s^2}{\sqrt{1-s^{p+1}}}ds 
= \frac{2}{\sqrt{k}}L_{p,2}\alpha^{(5-p)/2}.
\nonumber
\\
\frac12 &=& \int_0^{1/2} \frac{u_\lambda'(x)}
{\sqrt{k(\alpha^{p+1}-u_\lambda(x)^{p+1})}}dx 
\\
&=& \frac{1}{\sqrt{k}}\int_0^\alpha \frac{1}{\sqrt{\alpha^{p+1}-\theta^{p+1}}}d\theta
\nonumber
\\
&=& \frac{1}{\sqrt{k}}\alpha^{(1-p)/2}\int_0^1 \frac{1}{\sqrt{1-s^{p+1}}}ds
\nonumber
\\
&=& \frac{1}{\sqrt{k}}L_{p,0}\alpha^{(1-p)/2}.
\nonumber
\end{eqnarray}
By this, we have 
\begin{eqnarray}
\sqrt{k} = 2L_{p,0}\alpha^{(1-p)/2}.
\end{eqnarray}
By this, (3.7) and (3.8), we have 
\begin{eqnarray}
\Vert u_\lambda'\Vert_2^2 = 4L_{p,0}M_{p,2}\alpha^2, \quad \Vert u_\lambda\Vert_2^2 &=& \frac{L_{p,2}}{L_{p,0}}\alpha^2, \quad 
\Vert u_\lambda'\Vert_2^2 = \frac{4L_{p,0}^2M_{p,2}}{L_{p,2}}\Vert u_\lambda\Vert_2^2.
\end{eqnarray}
Thus the proof is complete. \qed

\vspace{0.2cm}

\noindent
{\bf Lemma 3.2.} {\it Assume that $u_\lambda$ satisfies (1.7). Then 
\begin{eqnarray}
\Vert u_\lambda'\Vert_2^2 =  \frac{4L_{p,0}^2M_{p,2}}{L_{p,2}}\Vert u_\lambda\Vert_2^2.
\end{eqnarray}
}

We obtain Lemma 3.2 by the same argument as that to prove Lemma 3.1. Therefore, we omit the proof. 

\vspace{0.2cm}

\noindent
{\bf Proof of Theorem 1.4.} Let $u_\lambda$ be a solution to (1.7). We put 
\begin{eqnarray}
d_0 := a\frac{4aL_{p,0}^2M_{p,2}}{L_{p,2}} + b.
\end{eqnarray}
Then by Lemmas 3.1 and 3.2, we have 
\begin{eqnarray}
d_0\Vert u_\lambda\Vert_2^2 + 1 = a\Vert u_\lambda'\Vert_2^2 + b\Vert u_\lambda\Vert_2^2 +1.
\end{eqnarray}
Therefore, the solution $u_\lambda$ of (1.7) is also the solution of (1.1), since (3.14) holds. 
Therefore, we are able to 
apply the argument in Section 2 and obtain Theorem 1.4. Thus the proof is complete. \qed

\section{Appendix}

 Let $p > 1$. We show (1.5) and 
 (1.6), which was proved in [16], for completeness. We apply the time map argument 
 to (1.3) (cf. [12]). Since (1.3) is autonomous, as (2.1)--(2.3), we have 
\begin{eqnarray}
W_p(x) &=& W_p(1-x), \quad 0 \le x \le \frac12,
\\
\xi& :=& \Vert W_p\Vert_\infty = \max_{0\le x \le 1}W_p(x) = W_p\left(\frac12\right),
\\
W_p'(x) &>& 0, \quad 0 \le x < \frac12.
\end{eqnarray}
By (1.3), for $0 \le x \le 1$, we have 
\begin{eqnarray}
\{W_p''(x) + W_p(x)^p\}W_p'(x) = 0.
\end{eqnarray}
By this and (4.2), we have 
\begin{eqnarray}
\frac12W_p'(x)^2 + \frac{1}{p+1}W_p(x)^{p+1} = \mbox{constant} = \frac{1}{p+1}
W_p\left(\frac12\right)^{p+1} 
= \frac{1}{p+1}\xi^{p+1}.
\end{eqnarray}
By this and (4.3), for $0 \le x \le 1/2$, we have, using $\theta = \xi s$,
\begin{eqnarray}
W_p'(x) = \sqrt{\frac{2}{p+1}(\xi^{p+1} - W_p(x)^{p+1})}.
\end{eqnarray}
By (4.1) and (4.6), we have 
\begin{eqnarray}
\frac12 &=& \int_0^{1/2} 1 dx = \int_0^{1/2} 
\frac{W_p'(x)}{\sqrt{\frac{2}{p+1}(\xi^{p+1} - W_p(x)^{p+1})}}dx 
\\
&=& \sqrt{\frac{p+1}{2}}\int_0^\xi \frac{1}{\sqrt{\xi^{p+1}-\theta^{p+1}}}d\theta
\nonumber
\\
&=& \sqrt{\frac{p+1}{2}}\xi^{(1-p)/2}\int_0^1 \frac{1}{\sqrt{1-s^{p+1}}}ds 
\nonumber
\\
&=& \sqrt{\frac{p+1}{2}}\xi^{(1-p)/2}L_{p,0}.
\nonumber
\end{eqnarray}
By this, we have 
\begin{eqnarray}
\xi = (2(p+1))^{1/(p-1)}L_{p,0}^{2/(p-1)}.
\end{eqnarray} 
This implies (1.6). We next show (1.5). By (4.1), (4.2), (4.6) and (4.8), we have 
\begin{eqnarray}
\Vert W_p'\Vert_m^m &=& 2\int_0^{1/2} W_p'(x)^{m-1}W_p'(x)dx 
\\
&=& 2\left(\frac{2}{p+1}\right)^{(m-1)/2}\int_0^{1/2} 
\left(\xi^{p+1} - W_p(x)^{p+1}\right)^{(m-1)/2}
W_p'(x)dx 
\nonumber
\\
&=& 2^{(m+1)/2}(p+1)^{-(m-1)/2}\int_0^\xi (\xi^{p+1} - \theta^{p+1})^{(m-1)/2}d\theta
\nonumber
\\
&=& 2^{(m+1)/2}(p+1)^{-(m-1)/2}\xi^{(m-1)(p+1)/2 + 1}\int_0^1 (1-s^{p+1})^{(m-1)/2}ds
\nonumber
\\
&=&  2^{mp/(p-1)}(p+1)^{m/(p-1)}L_{p,0}^{(mp+m-p+1)/(p-1)}M_{p,m}.
\nonumber
\end{eqnarray}
This implies (1.5). Thus the proof is complete. \qed


\begin{thebibliography}{20}
\labelsep=1em\relax

\bibitem[1]{1} C. O. Alves, F. J. S. A. Corr\'ea and T. F. Ma,  Positive solutions for a quasilinear elliptic equation of Kirchhoff type. Comput. Math. Appl. 49 (2005), 85--93.
\bibitem[2]{2} B. Cheng, {New existence and multiplicity of nontrivial solutions for nonlocal elliptic Kirchhoff type problems}, J. Math. Anal. Appl. {394} (2012), no. 2,  488--495.
\bibitem[3]{3} F. J. S. A. Corr\^{e}a, On positive solutions of nonlocal and nonvariational elliptic problems, Nonlinear Anal. 59 (2004), 1147--1155.
\bibitem[4]{4} F. J. S. A. Corr\^{e}a, D. C. de Morais Filho, On a class of nonlocal elliptic problems via Galerkin method. J. Math. Anal. Appl. 310 (2005), no. 1, 177--187.
\bibitem[5]{5} B. Gidas, W. M. Ni and L. Nirenberg, Symmetry and related properties via the maximum principle. Comm. Math. Phys. 68 (1979), 209--243.
\bibitem[6]{6} C. S. Goodrich, A one-dimensional Kirchhoff equation with generalized convolution coefficients. J. Fixed Point Theory Appl. 23 (2021), no. 4, Paper No. 73, 23 pp.
\bibitem[7]{7} C. S. Goodrich, A topological approach to nonlocal elliptic partial differential equations on an annulus. Math. Nachr. 294 (2021), 286--309.
\bibitem[8]{8} C. S. Goodrich, Differential equations with multiple sign changing 
convolution coefficients. Internat. J. Math. 32 (2021), no. 8, Paper No. 2150057, 28 pp.
\bibitem[9]{9}  C.S. Goodrich, An analysis of nonlocal difference equations with finite convolution coefficients. J. Fixed Point Theory Appl. 24 (2022), no. 1, Paper No. 1, 19 pp. 
\bibitem[10]{10} 
Z. Liang, F. Li and J. Shi, {Positive solutions of Kirchhoff-type non-local elliptic equation: a bifurcation approach}, Proc. Roy. Soc. Edinburgh Sect. A {147} (2017), no. 4, 875--894.
\bibitem[11]{11} F. Liu, H. Luo and G. Dai, {Global bifurcation and nodal solutions for 
homogeneous Kirchhoff type equations}, Electron. J. Differential Equations {2020} (29) (2020), 
pp. 1--13.
\bibitem[12]{12} {T. Laetsch}, {
The number of solutions of a nonlinear 
two point boundary value problem}, 
Indiana Univ. Math. J. {\bf 20}  1970/1971 1--13.
\bibitem[13]{13} Z. Liang, F. Li, J. Shi, Positive solutions to Kirchhoff type equations with nonlinearity having prescribed asymptotic behavior. Ann. Inst. H. Poincar\'e C Anal. Non Lin\'eaire 31 (2014), no. 1, 155--167. 
\bibitem[14]{14} O. M\'endez, On the eigenvalue problem for a class of Kirchhoff-type equations, J. Math. Anal. Appl. 494 (2021), no. 2, Paper No. 124671, 15 pp.
\bibitem[15]{15} T. Shibata, Bifurcation diagrams of one-dimensional Kirchhoff type equations, 
Adv. Nonlinear Anal. 12 (2023), 356--368.
\bibitem[16]{16} T. Shibata, Global and asymptotic behaviors of bifurcation curves of one-dimensional nonlocal elliptic equations, J. Math. Anal. Appl. 516 no.2 (2022), 126525. 
\bibitem [17]{17} T. Shibata,  Asymptotic behavior of solution curves of nonlocal 
one-dimensional elliptic equations, Bound. Value Probl. (2022), Paper No. 63.
\bibitem[18]{18} 
R. Sta\'nczy, Nonlocal elliptic equations, Nonlinear Anal. 47 (2001), 3579--3584.
\bibitem[19]{19} W. Wang and W. Tang, Bifurcation of positive solutions for a nonlocal problem, 
Mediterr. J. Math. 13 (2016), 3955--3964.


\end{thebibliography}
\end{document}